\documentclass{article}
\newcommand{\ww}{\omega}
\usepackage{amssymb}
\usepackage[english]{babel}
\usepackage{times}

\newtheorem{theorem}{Theorem}[section]

\newtheorem{e-proposition}[theorem]{Proposition}

\newtheorem{e-definition}[theorem]{Definition\rm}
\newtheorem{remark}{\it Remark\/}
\newtheorem{example}{\it Example\/}

\def\real{\mathbb{R}}
\def\C{\mathcal{C}}
\def\D{\mathcal{D}}
\def\E{\mathcal{E}}
\def\FF{\mathcal{F}}
\def\LL{\mathcal{L}}

\begin{document}
\title{\hrulefill\\
Cohomology of diffeological spaces and foliations \thanks{Partially supported by MEC Research Project MTM2004-05082 Spain}}
\author{E.~Mac\'{\i}as-Virg\'os; E.~Sanmart\'{\i}n-Carb\'on\\
\hrulefill}
\date{}
\maketitle

\medskip

\begin{abstract}
 Let $(M,\FF)$ be a foliated manifold. We study the relationship between  the basic cohomology $H_b(M,\FF)$ of the foliation and the
De Rham cohomology $H(\D_\FF)$ of the space of leaves  $M/\FF$
as a quotient diffeological space.
 We prove that for an arbitrary foliation there is a  morphism $H(\D_\FF)\to  H_b(M,\FF)$. It is an isomorphism when $\FF$ is a $Q$-foliation.  
\end{abstract}

\section{Introduction}
Diffeological spaces were introduced by J.-M.~Souriau  \cite{SOURIAU}  
(and in a slightly different form by K.-T. Chen 
\cite{CHEN}), 
as a generalization of the notion of manifold. Spaces of maps and quotient spaces naturally fit into this category, and many
definitions from differential geometry can be generalized to this setting. In particular, due to its contravariant nature,
the cohomology groups of a diffeological space can be defined in a  canonical way.
For a  complete study of the De Rham calculus in diffeological spaces see 
P.~Iglesias' book 
\cite{IGLESIAS}. 

From the point of view of foliation theory, the most interesting example is the space of leaves $M/\FF$ of a foliated manifold, endowed with the quotient diffeology $\D_\FF$. 
It is then a natural problem to compare its cohomology $H(\D_\FF)$ with the basic cohomology $H_b(M,\FF)$ of the foliation
\cite{REINHART}.

In this note we prove that there is an morphism $B\colon H(\D_\FF) \to H_b(M,\FF)$. Also we prove that $B$ is an isomorphism when $\FF$ is a $Q$-foliation (in the sense of R.~Barre \cite{BARRE}). This class of foliations includes  riemannian foliations without holonomy on a compact manifold;  totally geodesic foliations; and foliations having C.~Godbillon's homotopy extension property \cite{MEIGNIEZ}.

We thank F.~Alcalde for several useful remarks.
\section{Diffeological spaces}
 Let $X$ be a set.
Each map $\alpha \colon U \subset \real^n \to X$ defined on an open subset $U$ of some euclidean space $\real^n$, $n
\geq 0$, will be called a {\em plot} on $X$. A {\em diffeology of class $\C^\infty$} on $X$ is any family $\D$ of plots satisfying the following axioms:
\begin{enumerate}
\item All constant plots are in $\D$. 
\item If $\alpha \in \D$ is defined on $U\subset \real^n$, and $h\colon V \subset \real^m \to U$ is
any $\C^\infty$ map, then $\alpha \circ h \in \D$. 
\item Suppose that $\alpha \colon U \subset \real^n \to X$ is a plot and that each $t \in
U$ has a neighbourhood $V=V_t \subset U$ such that $\alpha\mid_{V}\in \D$. Then $\alpha \in \D$.
\end{enumerate}

A set $X$ endowed with a diffeology $\D$ is called a {\em diffeological space}. 

Let $(X, \D)$, $(Y , \E)$ be diffeological spaces. A map
$f\colon X \to Y $ is $\D$-{\em differentiable} if $f \circ \alpha \in \E$ for all $\alpha \in \D$.

\begin{example}[Finite dimensional manifolds]{\rm 
For instance, on a $\C^\infty$ manifold $M$ we can define the $\C^\infty$-{\em diffeology} $\D$ which consists on all the $\C^\infty$ plots $\alpha\colon U
\subset\real^n \to M$, $n\geq 0$. If $M$ and $N$ are $\C^\infty$ manifolds 
endowed with the  $\C^\infty$-diffeology, 
then a map $f\colon M \to N$ is $\D$-differentiable if and only if it is a map of class $\C^\infty$.
}
\end{example}

\begin{example}[Foliations as diffeological spaces]\label{LEVANTA}{\rm 
Let $(M,\FF)$ be a foliated manifold. The space of leaves $M/\FF$ (where two points on $M$ are identified iff they are on the same leaf of $\FF$) can
be endowed with the so-called {\em quotient diffeology} $\D_\FF$: it is the  least diffeology such that the canonical projection   $\pi \colon M \to
M/\FF$ is $\D$-diferenciable, when we consider on $M$ the $\C^\infty$-diffeology. 

Explicitly, 
a plot $\alpha \colon U \subset \real^n \to M/\FF$ belongs to $\D_\FF$ iff $\forall t \in U$ $\exists V_t \subset \real^n$, neighbourhood of $t$, such that
$ \alpha\mid_{V_t}= \pi \circ \beta_t$
 for some $ \C^\infty$-map $\beta_t \colon V_t \to M$.
 }
 \end{example}
%
\section{Cohomology of diffeological spaces}
Let $(X, \D)$ be a diffeological space. A $\D$-{\em differential form} 
of degree $r \geq 0$ on $X$ is a family ${\ww}=\{{\ww}_\alpha\}_{\alpha \in
\D}$ indexed by the plots $\alpha \colon U \subset \real^n \to X$, $\alpha \in \D$, where each ${\ww}_\alpha \in \Omega^r(U)$ is an $r$-form
in the domain $U$ of $\alpha$. We ask  the family ${\ww}$ to be compatible with the $\C^\infty$ maps $h\colon V\subset \real^m \to U$, in the sense
that
${\ww}_{\alpha \circ h}=h^*{\ww}_\alpha$.

Let us denote $\Omega^r(\D)$ the space of $\D$-differential $r$-forms.
The {\em exterior differential} $d \colon \Omega^r(\D) \to \Omega^{r+1}(\D)$ is given by $d{\ww} = \{d{\ww}_\alpha\}_{\alpha \in \D}$ if ${\ww} =
\{{\ww}_\alpha\}_{\alpha \in \D}.$ It is well defined.
Moreover it verifies $d\circ d =0$. 
The cohomology of the complex $(\Omega(\D),
d)$ will be called the {\em De Rham cohomology}  $H(\D)$ of the diffeological space $(X,\D)$.

Let $(X, \D)$, $(Y , \E)$ be diffeological spaces,   $f \colon X \to Y $ a $\D$-differentiable map. The pull-back
$f^*\colon \Omega^r(\E) \to \Omega^r(\D)$ given by
$ (f^*{\ww})_\alpha ={\ww}_{f\circ \alpha}$,  $\alpha\in\D$,
is well defined.
Moreover $f^*(d{\ww}) = d(f^*{\ww}), \forall {\ww} \in \Omega(\E)$, hence $f$ induces a map $f^*\colon
H(\E) \to H(\D)$ in cohomology.\\


\begin{theorem} \label{VARIEDAD}
Let $M$ be a   differentiable manifold of finite dimension $m < +\infty$, endowed with the  $\C^\infty$-diffeology $\D$. 
The   map
$F \colon \Omega(M) \to \Omega(\D)$ given by $F({\ww}) =\{\alpha^*{\ww}\}_{\alpha \in \D}$ induces an isomorphism
$H_{DR}(M) \cong H(\D)$ between the usual De Rham cohomology of $M$ and the cohomology of the diffeological space $(M,\D)$. 
\end{theorem}
{\bf Proof:\ } 
Clearly, the map
$F$ is well defined and
commutes with the exterior differential.
We define an inverse map 
 $G\colon \Omega^r(\D) \to \Omega^r(M)$ as follows.
 Let us choose an atlas  $\{\varphi_i
\colon V_i \subset M   \to  U_i \subset \real^m \}_{i \in I}$ on $M$. If $\theta \in \Omega^r(\D)$,  let  $G(\theta)$ be the only $r$-form on $M$ such that
$$G(\theta)\mid_{V_i} =  (\varphi_i)^{*}\theta_{\varphi_i^{-1}}, \quad \forall i \in I.$$
We left the reader to check that
the form $G(\theta) $ is smooth, 
well defined,
and commutes with the differential.

Finally, let us prove that $F$ and $G$ are inverse maps. If $\ww \in \Omega^r(M)$, then $G(F\omega)\mid_{V_i}={\ww}\mid_{V_i}$, for all $i\in I$.
On the other hand, let $\theta \in \Omega^r(\D)$. If 
$\alpha \in \D$ is a plot with domain  $U\subset \real^n$,  for each $t\in U$ we take a chart $\varphi_i
\colon V_i \subset M   \to  U_i \subset \real^m$  such that  $\alpha(t)\in V_i$. 
 Then we have, when restricted to $\alpha^{-1}(V_i)\subset U$, that 
 $$F(G(\theta))_\alpha=
 \alpha^*G(\theta)
 =\alpha^*\varphi_i^*\theta_{\varphi_i^{-1}} =
(\varphi_i \circ \alpha)^*\theta_{\varphi_i^{-1}}= \theta_{\varphi_i^{-1}\circ \varphi_i \circ
\alpha}=\theta_{\alpha}.$$ 

\section{Basic cohomology of foliations}
Let $(M,\FF)$ be a foliated manifold. A differential form $\omega \in \Omega^r(M)$ is {\em basic} if $i_X\omega=0$, $i_Xd\omega=0$, for any vector
field $X$ tangent to the foliation. In other words $\omega$ is horizontal (it vanishes in the tangent directions) and it is preserved by any tangent flow (the Lie derivative $\LL_X\omega = 0$). The {\em basic cohomology} $H_b(M,\FF)$ of $(M,\FF)$ is the cohomology of the complex $\Omega_b (M,\FF)$ of 
basic forms, endowed with the usual exterior differential $d$ on $M$.

\begin{theorem} \label{MORFISMO} There is a  morphism $H(\D_\FF) \to H_b(M,\FF)$ between the  cohomology  of the leaf space $M/\FF$, endowed with the quotient diffeology $\D_\FF$, and the basic cohomology of the foliation.   \end{theorem}

{\bf Proof:\ }
Accordingly to the proof of Theorem \ref{VARIEDAD}, we have the morphism
$B=G\circ \pi^*\colon \Omega(\D_\FF) \to \Omega(\D)\cong\Omega(M)$,
where $\pi^*$ is induced by the map $\pi\colon M \to M/\FF$.

{We must prove that if $\theta \in \Omega^r(\D_\FF)$, then $B(\theta)$ is a basic form.}
Take an atlas  $\{\varphi_i \colon V_i \subset M \to U_i \times W_i \subset \real^p \times \real^q  \}_{i\in I}$ adapted to the foliation  ($p=\dim\FF$).
For each $i\in I$, we have the following commutative diagram
$$\begin{array}{ccccc}
U_i \times W_i &\stackrel{\varphi_i^{-1}}{\to}& V_i& \stackrel{}{\subset} & M\\
p_2 \downarrow & &\downarrow \pi_i& &\downarrow \pi\\
W_i &\stackrel{\psi_i}{\to}& T_i &\stackrel{\tau_i}{\to}& M/\FF
\end{array}$$ 
\noindent where 
$\pi_i \colon V_i \to T_i$ is a submersion onto a local transverse manifold; $\psi_i \colon W_i \to T_i$ is the induced transverse chart;
and  $\tau_i \colon  T_i \to M/\FF$ is given by $ \tau_i(\overline{x})
= \pi(x)$ if  $\pi_i(x)=\overline{x}$. 
Then
$$B(\theta)\mid_{V_i} = \varphi_i^{*}\theta_{\pi  \circ \varphi_i^{-1}} = \varphi_i^{*}\theta_{\tau_i \circ \psi_i\circ p_2}
=(p_2 \circ \varphi_i)^*\theta_{\tau_i \circ \psi_i}.$$

Let $X$ be a vector field tangent to $\FF$. For any $x \in M$ let 
$i \in I$ be such that $x \in V_i$. 
Then
\begin{eqnarray*}
(i_XB(\theta))_x(v _1, \dots,v _{r-1}) &= &B(\theta)_x(X_x, v _1, \dots, v _{r-1})\\
&=& ((p_2 \circ \varphi_i)^*\theta_{\tau_i \circ \psi_i})_x(X_x, v _1,
\dots, v _{r-1})\\ &=& 0,
\end{eqnarray*}
because  $(p_2 \circ \varphi_i)_{*x}(X_x) = 0$.
A similar argument shows that $i_X dB(\theta) =0$.\\

\begin{remark} The morphism $\pi^*\colon \Omega(\D_\FF) \to \Omega(\D)$ is injective. \end{remark}
\section{$Q$-foliations}
A foliation $\FF$ in the manifold $M$ is called a {\em $Q$-foliation} when it has the following {\em $C^\infty$ tubular property} \cite{MEIGNIEZ}: for any manifold $X$ with a given base point $x$; any pair of maps $\beta,\gamma \colon X \to M$ such that $\pi\circ\beta =\pi\circ\gamma$; and any tangential path
$\mu \colon I \to M$ such that $\mu(0)=\beta(x)$, $\mu(1)=\gamma(x)$, there exists a neighbourhood $U\subset X$ of $x$ and a differentiable homotopy $H\colon U\times I \to M$ verifying:
\begin{itemize}
\item
$H(y,0)=\beta(y)$, $H(y,1)=\gamma(y)$, $\forall y\in U$;
\item
$H(x,t)=\mu(t)$, $\forall t\in I$;
\item
for any $y\in U$ the path $H(y,t)$ is tangential to $\FF$.
\end{itemize}

\begin{theorem} Let $\FF$ be a $Q$-foliation. Then the De Rham cohomology $H(\D_\FF)$ of the quotient diffeological space $M/\FF$ is isomorphic to
 the basic cohomology $H_b(M,\FF)$.
\end{theorem}
{\bf Proof:\ }
We follow the notation of Theorem \ref{MORFISMO}. Let $\omega\in \Omega^r_b(M,\FF)$ be a  basic $r$-form. In order to obtain $\theta\in\Omega(\D_\FF)$ such that $B(\theta)=\omega$ we have to define $\theta_\alpha$ for any plot $\alpha \colon U\subset \real^n \to M/\FF$ in $\D_\FF$.
We know  that $\forall x\in U$  $\exists V\subset U$ and $\exists\beta \in \D$ (both depend on $x$) such that  $\alpha\mid_{V}=\pi \circ
\beta$. Then we define
$$\theta_\alpha\mid_{V} = \theta_{\pi \circ \beta} = \beta^*{\ww}.$$
Clearly this form is smooth  and satisfies the compatibility condition.

Let us prove that it is  well defined.  Suppose that $V\cap W\neq \emptyset$ and that we take another lift $\gamma$ defined on $W$, that is $\alpha\mid_{W}=\pi\circ\gamma$.
Let $x\in V\cap W$.   
Choose an arbitrary tangential path $\mu$ in $M$ between $\beta(x)$ and $\gamma(x)$, and let $H\colon U\times I \to M$ be a homotopy  as in the definition of $Q$-manifold. Then
$\beta^*\omega =H^*_0\omega$ while $\gamma^*\omega=
H^*_1\omega$, so we just have to check that the function
$t\mapsto (H^*_t\omega)_x$
is constant. 

We shall prove that its derivative is null at $t=0$, the same argument being valid for any $t=t_0\in I$. 

The proof is inspired in the classical definition of Lie derivative  \cite[page 70]{WARNER}. Let $p=\mu(0)\in M$. Let $v=\mu^\prime(0)\in T_pM$ be the vector tangent to the curve
$\mu(t)$ at $t=0$. Let us denote 
$$L_v\omega=\lim_{t\to 0}{{1\over t}((H^*_t\omega)_x-(H^*_0\omega)_x)}$$ the derivative of the function above. In order to prove that $L_v\omega=0$
we shall proceed in several steps.
\begin{enumerate}
\item
If $r=0$, that is if $\omega$ is a basic function $f$, then
$L_vf=(d / dt)_{\mid t=0}(f\circ H_t)(x)= f_{*p}(v)=0$
because the vector $v$ is tangent to the foliation.
\item
If $\omega=df$ is  the differential of a basic function $f$, 
let $w\in T_x(V\cap W)$ be a tangent vector. It is defined by some curve  $\xi(s)$, that is $\xi(0)=x$ and
$\xi^\prime(0)=w$. Let us define the function $F(t,s) =fH_t\xi(s)=fH(t,\xi(s))$.
Then 
\begin{equation}\label{AQUI}
(L_v df)(w)=(d/dt)_{\mid t=0}(H_t^*df)_x(w) =({\partial^2 F / \partial_t\partial_s})_{\mid(0,0)}
\end{equation}
because
$$
(dH_t^*f)_x(w)= (f\circ H_t)_{*x}(w)   
=(d/ds)_{\mid s=0}(fH_t\xi(s))
.$$
For each  $s$, let  $v_s$ be the vector tangent  to the curve  $H(\xi(s),t)$ at $t=0$, which is tangent to the foliation.
Then (\ref{AQUI}) equals $({\partial^2 F / \partial_s\partial_t})_{\mid(0,0)}$, which is zero, because
$$({\partial F / \partial_t})_{\mid(0,s)}=
(d/dt)_{\mid t=0}(fH_t\xi(s))=  f_{*H(\xi(s),0)}(v_s)=0.$$
\item
 For any two forms $\omega$, $\mu$ we have $L_v (\omega\wedge \mu) =L_v\omega \wedge (H_0^*)_x\mu + (H_0^*)_x\omega \wedge L_v\mu$. In fact
\begin{eqnarray*}
&&H^*_t(\omega\wedge \mu) - H^*_0(\omega\wedge \mu)\\
&= &
 H^*_t\omega \wedge H^*_t\mu -  H^*_t\omega\wedge H^*_0\mu + H^*_t\omega\wedge H^*_0\mu
 - H^*_0\omega\wedge H^*_0\mu\\
 &=&  H^*_t\omega \wedge (H^*_t\mu - H^*_0\mu ) + (H^*_t\omega- H^*_0\omega)\wedge H^*_0\mu.
\end{eqnarray*}
\item
Finally, any basic form can be locally written (by means of an adapted chart $\varphi=(x^i,y^j)$) as a sum of products
$f(y^1,\dots,y^q)dy^{j_1}\wedge\cdots\wedge dy^{j_r}$ of basic functions and differentials of basic functions. Hence $L_v\omega=0$.
\end{enumerate}


\bigskip

{\small 
E.~Mac\'{\i}as-Virg\'os \\
{Department of Geometry and Topology, University of Santiago de Compostela,\\ Institute of Mathematics,
 Avda. Lope de Marzoa s/n. Campus Sur. 15782-Santiago de Compostela, Spain}\\
{\tt macias@zmat.usc.es}\\
{\tt http://www.usc.es/imat/quique}\\

\and E.~Sanmart\'{\i}n-Carb\'on\\
{Department of Mathematics, University of Vigo,\\
 F. CC.EE., R\'ua Leonardo da Vinci, Campus Lagoas-Marcosende. 36310-Vigo, Spain}\\
{\tt esanmart@uvigo.es}
}

\end{document}